\documentclass[12pt]{article}
\usepackage{graphicx}
\usepackage{latexsym}
\usepackage{amsfonts,amssymb,amsmath}
\newtheorem{thm}{Theorem}[section]
\newtheorem{cor}[thm]{Corollary}
\newtheorem{lem}[thm]{Lemma}
\newtheorem{pro}[thm]{Proposition}
\newtheorem{defn}[thm]{Definition}
\title{Two-letter group codes that preserve aperiodicity
        of inverse finite automata}

\author{ Jean-Camille Birget, \ Stuart W.\  Margolis
        \thanks{ \ Both authors were supported in part by NSF grant
                   DMS-9970471. The first author was also supported in 
                   part by NSF grant CCR-0310793. 
 The second author acknowledges the support of the Excellency Center,
``Group Theoretic Methods for the Study of Algebraic Varieties''
of the Israeli Science Foundation  }
       }

\date{}
\begin{document}
\maketitle


\begin{abstract}
We construct group codes over two letters (i.e., bases of
subgroups of a two-generated free group) with special properties.
Such group codes can be used for reducing algorithmic problems
over large alphabets to algorithmic problems over a two-letter
alphabet. Our group codes preserve aperiodicity of inverse finite
automata. As an application we show that the following problems
are {\sc PSpace}-complete for {\em two-letter} alphabets (this was
previously known for large enough finite alphabets): The
intersection-emptiness problem for inverse finite automata, the
aperiodicity problem for inverse finite automata, and the
closure-under-radical problem for finitely generated subgroups of
a free group. The membership problem for 3-generated inverse
monoids is {\sc PSpace}-complete.
\end{abstract}


\section{Introduction}

Codes and coding theory are a well-known and important subject. In
its most general form, a code over an alphabet $A$ is defined to
be a subset $C$ of $A^*$  such that any concatenation of elements
of $C$ can be uniquely factored, or ``decoded'', into a sequence
of elements of $C$. Equivalently, the submonoid $\langle C
\rangle$ of $A^*$ generated by $C$ is free with base $C$, i.e.,
$\langle C \rangle$  is isomorphic to the free monoid $C^*$.
As a reference see \cite{BP}. Some notation: $A^*$ denotes the
free monoid over $A$, i.e., the set of all finite sequences
(``words'') of elements of $A$, including the empty word. $A^+$
denotes the free semigroup over $A$, i.e., the set of all
non-empty finite sequences over $A$.

For groups one can use the same definition of a code, replacing
``free monoid'' by ``free group''. In the literature such a code
is called a {\it base of a free group}. We'll call it {\it group
code} because we will use it in the spirit of information coding.
A precise definition of a group code appears below. First we need
some notation: The free group over a generating set $A$ is denoted
by FG$(A)$. We use a copy $A^{-1} = \{a^{-1} : a \in A\}$ of $A$,
disjoint from $A$, to denote the inverses of the generators. We
denote $A \cup A^{-1}$ by $A^{\pm 1}$. For $w = a_1 \ldots a_n$
with $a_1, \ldots, a_n \in A^{\pm 1}$, the inverse of $w$ is
defined to be $w^{-1} = a_n^{-1} \ldots a_1^{-1}$, where
$(a^{-1})^{-1}$ is always replaced by $a$ for all $a \in A$. The
identity element of FG$(A)$ is the empty word, and is denoted by
1. The elements of FG$(A)$ are all the words over the alphabet
$A^{\pm 1}$ that are {\it reduced}, i.e., that contain no
subsegment of the form $a \, a^{-1}$ or $a^{-1}a$ (for any $a \in
A$). In general, for any word $w \in (A^{\pm 1})^*$ we define the
{\it reduction of} $w$ to be the word obtained by cancelling all
subsegments of the form $w \, w^{-1}$ (with $w \in (A^{\pm 1})^*$)
iteratively as much as possible, and we denote the resulting
reduced word by red$(w)$. For any word $w$ we denote its length by
$|w|$. See \cite{MKS, LyS, Cohen} for background on free groups.

Any function $f: A \to (B^{\pm 1})^*$ can be extended (uniquely)
to a group morphism $f^{(G)}: {\rm FG}(A) \to {\rm FG}(B)$ defined
by $f^{(G)}(a_1^{\varepsilon_1} \ldots a_n^{\varepsilon_n})$ $ =
{\rm red}(f(a_1)^{\varepsilon_1} \ldots $
$f(a_n)^{\varepsilon_n})$, where $\varepsilon_1, \ldots,
\varepsilon_n \in \{-1, 1\}$.

{\it Important convention:} Throughout this paper we will view the
free group ${\rm FG}(A)$ as a subset of the free monoid $(A^{\pm
1})^*$; namely, ${\rm FG}(A)$ consists of all the reduced words
over $A^{\pm 1}$. Of course, ${\rm FG}(A)$ is only a subset of
$(A^{\pm 1})^*$, not a submonoid.

\begin{defn}  \   
Let $\varphi: A \to (B^{\pm 1})^*$ be a map whose extension to a
free-group morphism $\varphi^{(G)}: {\rm FG}(A) \to {\rm FG}(B)$ \  
is {\em injective}. Then the image set $\varphi^{(G)}(A)$ $( \,  
\subset {\rm FG}(B) \subset (B^{\pm 1})^*)$ is called a {\bf group
code} over $B$, and the elements of $\varphi^{(G)}(A)$ are called
{\bf code words}. By our convention, ${\rm FG}(B)$ is a subset of
$(B^{\pm 1})^*$, and hence a group code is a set of words.

The injective map $\varphi^{(G)}|_A: A \to {\rm FG}(B)$ defined by
\ $a \mapsto {\rm red}(\varphi(a))$, i.e., the restriction of
$\varphi^{(G)}$ to $A$, is called a {\bf group encoding} of $A$
over $B$.
\end{defn}
The study of free groups and of bases of free groups (i.e., group
codes) has a long history \cite{MKS, LyS, Cohen}. In particular,
Nielsen showed in the 1920s that every finitely generated subgroup
of a free group is itself free and hence has a group code. A
little later in the 1920s Schreier extended Nielsen's result to
all subgroups of a free group. So, group codes can be finite or
infinite. We note the following however:

\begin{pro} \  
An infinite group code cannot be a regular language, but can be
deterministic context-free.
\end{pro}
{\bf Proof.} \ If an infinite regular group code existed we could
apply the Pumping Lemma, so the group code would contain all words
of the form $w_n = ux^nv$ ($n \in {\mathbb N}$), for some fixed
words $u, x, v$, with $x$ non-empty. But then the following
non-trivial relation would hold among code words: \ $w_2 \, 
w_1^{-1} w_2 = w_3$.

The example $\{ a^n b a^{-n} : n \geq 0 \}$ over the alphabet
$\{a,b\}^{\pm 1}$, shows that there are infinite group codes that
are deterministic context-free languages. The set $\{a^n b a^{-n}
: n \geq 0\}$ is a well-known Nielsen basis.  \ \ \ $\Box$

\medskip

We are interested in group codes over an alphabet of size 2. Just
as for the usual codes (over monoids), the main purpose of group
codes is to translate large alphabets into smaller alphabets. This
in turn can be used to show that some problems that are hard over
large alphabets are also hard over a 2-letter alphabet. We will
consider the fixed two-letter alphabet $\{a,b\}$ and the inverses
$a^{-1}, b^{-1}$ of these letters.

Subgroups of a free group are closely related to inverse monoids
and inverse finite automata \cite{MM2}. By definition, an {\bf
inverse finite automaton} is a structure ${\cal A} = (Q, X,
\delta, q_0, q_f)$ where, according to the standard notation in
\cite{HU}, $Q$ is the set of states, $q_0$ is the start state, and
$q_f$ is the accept state. For inverse automata, the input
alphabet is $X \cup X^{-1} = X^{\pm 1}$, although we only mention
$X$ explicitly; the designation ``inverse'' automatically provides
the inverse letters. The state-transition relation $\delta$ is a
partial function $\delta: Q \times X^{\pm 1} \to Q$, and is
required to have the following property: For each letter $x \in
X$, the partial function $\delta(\cdot, x): q \in Q \mapsto
\delta(q, x) \in Q$ is {\it injective}. Moreover, we require that
the partial function $\delta(\cdot, x^{-1})$ be the inverse of
$\delta(\cdot, x)$. We represent an inverse finite automaton by
its state-graph, in the same way as for ordinary finite automata
(see \cite{HU}), except that we omit the edges labeled by inverse
letters. More precisely, when $\delta(p, x) = q$ (with $p,q \in
Q$, $x \in X$) we draw an edge $p \, {\buildrel {x} \over
\longrightarrow} \, q$; we implicitly also have an edge $q \,
{\buildrel {x^{-1}} \over \longrightarrow} \, p$, but we don't
draw that edge. See e.g. \cite{BMMW} for more information on
inverse automata.

Let $\kappa : X^{\pm 1} \to (\{a,b\}^{\pm 1})^*$ be any group
encoding and let ${\cal A}$ be any inverse finite automaton ${\cal
A}$ with input alphabet $X$. We define the {\bf encoded inverse
finite automaton} $\kappa({\cal A})$, with input alphabet
$\{a,b\}$, by the following two-step construction:

\smallskip

\noindent (1) \ We replace every edge \ $p \, {\buildrel {x} \over
\longrightarrow} \, q$ \ of ${\cal A}$ (with $x \in X$) by a path
labeled by $\kappa(x)$; to do this we introduce $|\kappa(x)|-1$
new states and $|\kappa(x)|$ new edges. Implicitly, we now also
have the inverses of the new edges, thus obtaining a path from $q$
to $p$ labeled by $\kappa(x^{-1})$.  Let $\kappa({\cal A})_0$ be
the nondeterministic finite automaton obtained so far.

\smallskip

\noindent (2) \ Starting from $\kappa({\cal A})_0$ we apply the
{\em fold} operation as much as possible. This means that any two
edges (explicitly drawn or implicit) with a common beginning or
end vertex, and with identical label in $\{a,b\}^{\pm 1}$ are made
equal. For example, if
 \ $p \, {\buildrel {x^e} \over \longrightarrow} \, q_1$ \ and
 \ $p \, {\buildrel {x^e} \over \longrightarrow} \, q_2$ \ are present (with
$e \in \{-1,1\}$) then one folding step makes $q_1$ equal to
$q_2$, and the above two edges become equal. See e.g., \cite{Sta},
\cite{MM2}, \cite{BMMW} for more information on the very classical
fold operation. In particular, it is well known that maximal
folding produces a unique resulting automaton, which does not
depend on the folding sequence chosen. We denote this resulting
automaton by $\kappa({\cal A})$; it is an inverse automaton if
${\cal A}$ is an inverse automaton. We denote the transition
function of $\kappa({\cal A})$ by $\delta_{\kappa}$.

\smallskip

In general, for any automaton ${\cal M}$ we let $L_{\cal M}$
denote the language accepted by ${\cal M}$. For an inverse
automaton ${\cal A} = (Q, A, \delta, q_0, q_f)$ we consider the
language accepted $L_{\cal A} \subseteq (A^{\pm 1})^*$, as well as
the group language of ${\cal A}$, defined as follows:

\begin{defn} \  
 The {\bf group language} of a finite inverse automaton
${\cal A}$ with input alphabet $A$ consists of the {\em reduced}
words ($\in (A^{\pm 1})^*$) accepted by ${\cal A}$; in other
words, the group language of ${\cal A}$ is $L_{\cal A} \cap {\rm
FG}(A)$.
\end{defn}

\begin{lem} \label{Lem1.1} \  
For a finite inverse automaton ${\cal A}$ with input alphabet $A$
the group language \ $L_{\cal A} \cap {\rm FG}(A)$ =
 \ ${\rm red}(L_{\cal A})$.
\end{lem}
{\bf Proof.} This is Lemma 1.1 in \cite{BMMW}. \ \ \ $\Box$

\medskip

Note that by Benois' theorem \cite{Benois}, \cite{Berstel}, ${\rm
red}(L_{\cal A})$ is also accepted by a finite automaton with
alphabet $A^{\pm 1}$. But this automaton cannot be an inverse
automaton, except in trivial cases. Indeed, an inverse automaton
will always accept some non-reduced words (except when $L_{\cal
A}$ is empty or consists of just the empty word).

An {\em automaton with involution} over the alphabet $(A^{\pm 1})^*$ is
an automaton ${\cal A}$ such that for every edge $p \, {\buildrel
{x} \over \longrightarrow} \, q$ with $x \in (A^{\pm 1})^*$, of
${\cal A}$, $q \, {\buildrel {x^{-1}} \over \longrightarrow} \, p$
is also an edge of ${\cal A}$. We will always assume that all
automata over the alphabet $A^{\pm 1}$ are automata with
involution. Notice that an automaton with involution is
deterministic if and only if it is an inverse automaton.

Let ${\cal A}$ be any automaton with involution over the alphabet
$A^{\pm 1}$. The folded automaton ${\rho(\cal A)}$ is defined as
above by applying some maximal folding sequence to ${\cal A}$.
This determines an equivalence relation $\sim$ on the states of
${\cal A}$ by defining two states to be equivalent if they define
the same state of ${\rho(\cal A)}$, that is, if the two states are
folded onto one another. Recall that a Dyck word over $(A^{\pm
1})^*$ is a word that reduces to the identity word in ${\rm
FG}(A)$. The language of Dyck words is known to be the smallest
language containing the empty word and closed under concatenation
and the conjugation operation $w \mapsto awa^{-1}$, for all 
$a \in A^{\pm 1}$.

\begin{lem} \label{Dyck} \   
Let ${\cal A}$ be an automaton with involution over the alphabet
$(A^{\pm 1})$. Then states $p,q$ of ${\cal A}$ satisfy $p \sim q$
if and only if there is a Dyck word $w$ such that $w$ labels a
path from $p$ to $q$ in ${\cal A}$.
\end{lem}

{\bf Proof.} Assume that the reduced automaton ${\rho(\cal A)}$ is
obtained by a sequence of $m$ foldings. Let ${\cal A}_{i}$ be the
automaton obtained after $i$ foldings, $0 \leq i \leq m$. There is
a corresponding equivalence relation $\sim_i$ on the states of
${\cal A}$, and
 \ $\sim_{0} \ \subset \ \sim_{1} \ \subset \ \ldots \ $   
 $ \subset \ \sim_{m} \ =  \ \sim$.

We will prove by induction  that if  $i$ is the least integer
such that $p \sim_i q$, then there is a Dyck word $w$ that
 labels a path from $p$ to $q$ in ${\cal A}$. This is true if
$i=0$ since then the empty word labels a path from $p$ to itself.

Assume that if $r \sim_i s$ then there is a Dyck word labeling 
a path from $r$ to $s$ in ${\cal A}$; and assume that $p
\sim_{i+1} q$, but $p \not\sim_i q$. Since a folding identifies
exactly two states, the ($i+1$)st folding identifies the $\sim_i$
class of $p$ with that of $q$. Let $[r]_{\sim_i}$ denote the
$\sim_i$ equivalence class of a state $r$ of ${\cal A}$.

Thus there is a $\sim_i$ equivalence class, $X$, such that there
are edges of ${\cal A}_i$, $[p]_{\sim_i} \, {\buildrel {x} \over
\longrightarrow} \, X$ and $X \, {\buildrel {x} \over
\longleftarrow} \, [q]_{\sim_i}$ for some $x \in A^{\pm 1}$. It
is clear that every path in ${\cal A}_{i}$ lifts, by ``unfolding'',
 to a path of ${\cal A}$. 
Thus in ${\cal A}$ there are states $p',q'$ and states
$r,s \in X$ such that $p' \in [p]_{\sim_i}$, 
$q' \in [q]_{\sim_i}$ and    
$p' \, {\buildrel {x} \over \longrightarrow} \, r$ and $s \,
{\buildrel {x} \over \longleftarrow} \, q'$ in ${\cal A}$. Since
 \ $p \sim_{i} p' \, {\buildrel {x} \over \longrightarrow} \, $
$r \sim_{i} s \, {\buildrel {x} \over \longleftarrow} \, $
$q \sim_{i} q'$ \ we have, by induction, 
Dyck words $u,v,w$ that label paths from $p$ to $p'$, $q'$ to $q$
and $r$ to $s$ respectively in ${\cal A}$. Therefore the Dyck word
$uxwx^{-1}v$ labels a path from $p$ to $q$ in ${\cal A}$.

Conversely, a straightforward induction on the length of a Dyck
word $w$ shows that  if $w$ labels a path from a state $p$ to a state
$q$ of ${\cal A}$ then $p \sim q$. $\Box$

\begin{cor} Let ${\cal A}$ be an automaton with involution over the 
alphabet
$A^{\pm 1}$ and let $\rho({\cal A})$ be the reduced automaton of
${\cal A}$. Let $p,q$ be states of ${\cal A}$. If $w = a_{1}
\ldots a_{n}$, with $a_{i} \in A^{\pm 1}, 1 \leq i \leq n$, labels a
path from $[p]_{\sim}$ to $[q]_{\sim}$ in $\rho({\cal A})$, then
there are Dyck words $u_{0}, \ldots, u_{n}$ such that $u_{0}a_{1}
\ldots a_{n}u_{n}$ labels a path from $p$ to $q$ in ${\cal A}$. In
particular, ${\rm red}\big(L({\cal A})\big) = $
${\rm red}\big(L(\rho({\cal A}))\big)$.
\end{cor}

{\bf Proof.} There are states $p = p_{0}$, $p_{1}$, $ \ldots,$ 
$ p_{n} = q$ of ${\cal A}$ such that 
$[p_{i}]_{\sim} \, {\buildrel {a_{i+1}}
\over \longrightarrow} \, [p_{i+1}]_{\sim}$ are edges of
$\rho({\cal A})$. Since paths in $\rho({\cal A})$ lift to paths of
${\cal A}$, there are states $p'_{0}, p'_{1}, \ldots, p'_{n}$ of
${\cal A}$ such that $p_{i} \sim p'_{i}$ for $0 \leq i \leq n$, and
such that there are
edges $p'_{i} \, {\buildrel {a_{i+1}} \over \longrightarrow} \, 
p'_{i+1}$ of ${\cal A}$. By Lemma \ref{Dyck}, there are Dyck words
$u_{0}, \ldots, u_{n}$ such that $p_{i} \, {\buildrel {u_{i}}
\over \longrightarrow} \, p'_{i}$ and the first assertion of the
corollary follows.

It is clear that ${\rm red}\big(L({\cal A})\big) \subseteq {\rm
red}\big(L(\rho({\cal A}))\big)$ since paths in ${\cal A}$ fold to 
paths in $\rho({\cal A})$. 
The converse inclusion follows from the first
assertion of the corollary if we take $w$ to be a reduced word.
$\Box$

We record a special case of the above corollary that is of special
interest in this paper in the proposition below.
\smallskip

\begin{pro} \label{preservCode} \  
Let $\kappa: X \to (A^{\pm 1})^*$ be any group encoding, and let
$\kappa^{(M)}: (X^{\pm 1})^* \to (A^{\pm 1})^*$ be the
corresponding monoid morphism. Let ${\cal A}$ be an inverse finite
automaton with alphabet $X$ and let $L_{\cal A} \subseteq (X^{\pm
1})^*$ be the language it accepts. Then the group language of
$\kappa({\cal A})$ is ${\rm red}\big(\kappa^{(M)}(L_{\cal A})\big)$. 
In other words, \ ${\rm red}\big(L_{ \kappa({\cal A}) }\big) \ = \   
{\rm red}\big(\kappa^{(M)}(L_{\cal A})\big)$.
\end{pro}

\section{Aperiodicity preserving group codes}

Some standard definitions: A monoid $M$ is called {\em aperiodic}
iff $x^{n+1} = x^n$ for all $x \in M$, for some constant $n$
depending only on $M$.  A finite automaton ${\mathcal A}$ is
called aperiodic iff the syntactic monoid of ${\mathcal A}$ is
aperiodic.

Let $Y$ be a finite subset of {\rm FG}$(A)$, and let $H = \langle
Y \rangle$ be the subgroup of {\rm FG}$(A)$ generated by $Y$. Then
we can construct a finite inverse automaton ${\mathcal A}_H$ with
the following property: A reduced word $w \in {\rm FG}(A)$ belongs
to $H = \langle Y \rangle$ iff ${\mathcal A}_H$ accepts $w$. In
other words:  \ The group language $L({\mathcal A}_H) \cap {\rm
FG}(A)$ of ${\mathcal A}_H$ is $H$. A construction of ${\mathcal
A}_H$ goes as follows (see \cite{BMMW}, p.\ 251, for more
details): Consider cyclic graphs labeled by the elements of $Y$,
and glue these cycles together at their origins; if we now pick
this common origin as the start and accept state we obtain a
nondeterministic automaton. Next, we apply maximal folding. The
resulting finite inverse automaton is ${\mathcal A}_H$. One can
show that it only depends on $H$ (not on the originally given
generating set $Y$).

\begin{defn} \   
A subgroup $H$ of a group $G$ is {\bf closed under radical} (also
called ``radical-closed'', or ``pure'') iff for all $g \in G$ and
all $N > 0$ we have: \ $g^N \in H$ implies $g \in H$.

The radical of $H$ in $G$ is the set \ $\sqrt{H} = \{ g \in G : $
there exists $N > 0$ with
 $g^N \in H \}$.
\end{defn}

Closure under radical for subgroups of a free group is intimately
connected to aperiodicity of inverse automata:

\begin{lem} \label{Aperiodicity_Purity} \  
Let $Y$ be a finite subset of {\rm FG}$(A)$. The subgroup $H =
\langle Y \rangle$ of {\rm FG}$(A)$ generated by $Y$ is closed
under radical iff the finite inverse automaton ${\mathcal A}_H$ is
aperiodic.
\end{lem}
{\bf Proof.} This is Theorem 3.1 in \cite{BMMW}. \ \ \ $\Box$

\medskip

\begin{pro} {\bf (Transitivity of radical closure).} \   
Consider subgroups $K \leq H \leq G$ such that $K$ is
radical-closed in $H$ and $H$ is radical-closed in $G$. Then $K$
is radical-closed in $G$.
\end{pro}
{\bf Proof.} Suppose $x \in G$ is such that $x^n \in K$, for some
integer $n \geq 2$. Then $x^n \in H$, hence $x \in H$, by radical
closure of $H$ in $G$. So we have now $x \in H$ and $x^n \in K$.
This implies that $x \in K$, by radical closure of $K$ in $H$. 
 \ \ \ $\Box$

\begin{defn} \  
A group homomorphism $h: {\rm FG}(X) \to {\rm FG}(A)$ preserves
closure under radical iff for every subgroup $H$ of ${\rm FG}(X)$
we have: \ $H$ is closed under radical in ${\rm FG}(X)$ iff $h(H)$
is closed under radical in ${\rm FG}(A)$.

A group encoding $\varphi: X \to (A^{\pm 1})^*$ is said to
preserve closure under radical iff the group homomorphism
$\varphi^{(G)}: {\rm FG}(X) \to {\rm FG}(A)$ determined by
$\varphi$ preserves closure under radical.
\end{defn}

\begin{pro}  \label{preserv} \   
Let $f: {\rm FG}(X) \to {\rm FG}(A)$ be an injective morphism such
that the image group Im($f$) of $f$ is radical-closed in ${\rm
FG}(A)$. Then for all subgroups $H$ of ${\rm FG}(X)$ we have: \  
$H$ is radical-closed in ${\rm FG}(X)$ iff $f(H)$ is
radical-closed in ${\rm FG}(A)$. \ In other words:

A group encoding $\varphi$ preserves radical-closure iff
Im($\varphi$) (reduced in the free group) is radical-closed.
\end{pro}
{\bf Proof.} Suppose $f(H)$ is radical-closed in ${\rm FG}(A)$.
Then $f(H)$ is also radical-closed in Im($f$). Hence, since $f$ is
an isomorphism between the groups ${\rm FG}(X)$ and Im($f$), $H$
is radical-closed in ${\rm FG}(X)$.

Suppose $H$ is radical-closed in ${\rm FG}(X)$. Then $f(H)$ is
radical-closed in Im($f$), since $f$ is an isomorphism between
${\rm FG}(X)$ and Im($f$). Hence, since Im($f$) is radical-closed
in ${\rm FG}(A)$, transitivity of radical closure implies that
$f(H)$ is also radical-closed in ${\rm FG}(A)$.  \ \ \ $\Box$

\bigskip

\noindent {\bf Example:} {\it A family of finite aperiodic
two-letter group codes of all sizes}

\medskip

Consider \ $C_n = \{ a^i b a^{-i} : 0 \leq i \leq n-1\}$, over the
alphabet $\{a,b\}^{\pm 1}$. It is well known that this set has the
Nielsen property, hence it is a group code (compare with Ex.\ 3,
Sect. 3.2, p.\ 138 in \cite{MKS}). Moreover, the inverse automaton
${\mathcal A}$ given by the following transition table (with state
set $\{1,2, \ldots, n\}$, with 1 as both start and accept state)
satisfies:

\smallskip

red$(L_{\mathcal A}) = $ red($\langle C_n \rangle$),

\smallskip

\noindent where ``red'' refers to reduction in FG$(\{a,b\})$. In
other words, the free group red($\langle C_n \rangle$) is the
group language of ${\mathcal A}$.

\bigskip

\hspace{1in}
\begin{tabular}{l||l|l|l|l|l|}
    & 1 & 2 & \ldots & $n-1$ & $n$ \\ \hline \hline
$a$ & 2 & 3 & \ldots & $n$ & $-$ \\ \hline $b$ & 1 & 2 & \ldots &
$n-1$ & $n$ \\ \hline
\end{tabular}

\bigskip

The syntactic inverse monoid of ${\mathcal A}$ is generated by the
identity map, corresponding to the letter $b$, and the partial map
 \ $i \in \{1,2, \ldots, n-1\} \longmapsto i+1$ \ (undefined on $n$),
corresponding to the letter $a$. Since this is a one-generator
monoid with zero, satisfying $a^n = 0$, the monoid is aperiodic.

In summary we have:

\begin{pro} \label{mainProp} \   
For any alphabet $X = \{x_1, x_2, \ldots, x_n \}$ of size $n$, the
map $f: x_i \mapsto a^{i-1} b a^{-i+1}$ ($1 \leq i \leq n$) is a
group encoding into a two-generated free group that preserves
closure under radical.
\end{pro}

By combining the above lemmas and propositions we obtain:
\begin{cor} \label{codeAper} \   
Let $f$ be the group encoding defined in Proposition
\ref{mainProp}. Let $\{ w_1, \ldots, w_k\}$ be any finite set of
words $\subset (X^{\pm 1})^*$. Then the subgroup $\langle w_1,
\ldots, w_k \rangle$ of ${\rm FG}(X)$ is closed under radical iff
the subgroup $\langle f(w_1), \ldots, f(w_k) \rangle$ of ${\rm
FG}(\{a,b\})$ is closed under radical.
\end{cor}

\medskip

\noindent {\bf Application: Complexity of radical-closure
       and aperiodicity problems}

\medskip

Group encodings are log-space computable reductions from large
alphabets to small alphabets. This enables us to show that the
problems below  about inverse finite automata or about free groups
are {\sc PSpace}-complete over two-letter alphabets. Previously it
was known that they are {\sc PSpace}-complete over all large
enough finite alphabets (\cite{BMMW}, Theorem 6.13).

The {\em aperiodicity problem} takes as input a finite automaton
and asks whether the language accepted by this automaton is
aperiodic. S.~Cho and D.~Huynh \cite{ChoHuynh} showed that the
aperiodicity problem for general finite automata is {\sc
PSpace}-complete, and it was shown in \cite{BMMW} (Theorem 6.13)
that the problem remains {\sc PSpace}-complete for inverse finite
automata (over some finite alphabet).

The {\em radical-closure problem} for a free group FG$(X)$ takes
as input a list of words $w_1, \ldots, w_n$ $\in$ FG$(X)$, and
asks whether the subgroup $\langle w_1, \ldots, w_n \rangle$ of
FG$(X)$ generated by these words is closed under radical. It was
proved in \cite{BMMW} (Theorem 7.1) that this problem is {\sc
PSpace}-complete for some finite alphabet $X$. We can now
strengthen these results:

\medskip

\begin{thm} \label{PSpace_compl_Aper} \  
The radical-closure problem for a free group with {\em two
generators}, and the aperiodicity problem for inverse finite
automata over a {\em two-letter} alphabet, are {\sc
PSpace}-complete.
\end {thm}
{\bf Proof.} \ By Corollary \ref{codeAper}, the group encoding $f$
is a reduction of the radical-closure problem over any fixed
finite alphabet to the radical-closure problem over a two-letter
alphabet. It was shown in \cite{BMMW} (Theorem 3.6) that the
radical-closure problem and the aperiodicity of inverse finite
automata are polynomial-time reducible to each other; in this
reduction, the alphabets are preserved.

Finally, as we saw above, the radical-closure problem is {\sc
PSpace}-complete over some finite alphabet, and is in {\sc PSpace}
for all finite alphabets.
 \ \ \ $\Box$


\section{Other applications of group codes}

\bigskip

\noindent As we saw, a group encoding is a log-space computable
function from a possibly large alphabet problems to a possibly
small alphabet. This will enables us to show that the problems
below  about inverse finite automata or about free groups are {\sc
PSpace}-complete over a two- or three-letter alphabet.

The {\bf intersection-emptiness problem} for finite automata takes
as input a list of finite automata ${\cal A}_i$ ($i = 1, \ldots,
k$) where $k$ is part of the input, and asks whether the
intersection of the languages accepted by these automata is empty.
For general deterministic finite automata this problem was shown
to be {\sc PSpace}-complete by D.~Kozen \cite{Kozen}, and for
inverse finite automata  {\sc PSpace}-completeness was shown in
\cite{BMMW} (Proposition 5.3).

\medskip

\begin{thm} \label{PSpace_compl} \    
The intersection-emptiness problem for inverse finite automata
over a fixed {\em two-letter} alphabet is {\sc PSpace}-complete.
\end {thm}
{\bf Proof.} \ Let ${\cal A}_1, \ldots, {\cal A}_n$ be inverse
finite automata with alphabet $A$ and let $L_1, \ldots, L_n
\subseteq (A^{\pm 1})^*$ be the respective languages that they
accept. Let $f: A \to (B^{\pm 1})^*$ be any group encoding with
$|B| =2$, and let $L'_1, \ldots, L'_n \subseteq (B^{\pm 1})^*$ be
the languages accepted by the inverse finite automata $f({\cal
A}_1), \ldots, f({\cal A}_n)$ respectively.

We claim that $L_1 \cap \ldots \cap L_n = \varnothing$ iff $L'_1
\cap \ldots \cap L'_n = \varnothing$, which shows that $f$ reduces
the intersection emptiness problem of inverse automata over the
alphabet $A$ to the intersection emptiness problem of inverse
automata over the alphabet $B$.

If $L_1 \cap \ldots \cap L_n \neq \varnothing$ consider $w \in L_1
\cap \ldots \cap L_n$. By Lemma \ref{Lem1.1} we can assume that
$w$ is reduced. Then, by Prop.\ \ref{preservCode}, \ ${\rm
red}(f(w)) \in L'_1 \cap \ldots \cap L'_n$; hence $L'_1 \cap
\ldots \cap L'_n \neq \varnothing$.

Conversely, if $y \in L'_1 \cap \ldots \cap L'_n$ $(\neq
\varnothing)$ we can again assume by Lemma \ref{Lem1.1} that $y$
is reduced. Then by Prop.\ \ref{preservCode},
 \ $y \in {\rm red}(f(L_1)) \cap \ldots \cap {\rm red}(f(L_n))$. Since the
function \ $F = {\rm red}(f(.)): {\rm FG}(A) \to {\rm FG}(B)$ is
injective (by definition of a group code), it has an inverse
function $F^{-1}$ and we have $F^{-1}(y) \in L_1 \cap \ldots \cap
L_n$.  So, $L_1 \cap \ldots \cap L_n \neq \varnothing$.

Finally, as we saw above, the intersection-emptiness problem is
{\sc PSpace}-complete over some finite alphabet. So the reduction
makes the encoded problems {\sc PSpace}-complete over a two-letter
alphabet.
 \ \ \ $\Box$

\bigskip

The {\bf membership problem} for finite inverse monoids is defined
as follows: The input is a finite list of injective partial maps
$f_0, f_1, \ldots, f_m$ on a finite set $\{1, \ldots, n\}$. Each
$f_i$ is described by a function table that bijectively maps a
subset of $\{1, \ldots, n\}$ to a subset of $\{1, \ldots, n\}$;
entries in the table where $f_i$ is not defined are blank. The
question is whether $f_0$ can be written as a composition of some
of the $f_i$ and $f_i^{-1}$ ($1 \leq i \leq m$); more rigorously,
the question is  whether $f_0$ belongs to the inverse monoid
generated by $\{f_1, \ldots, f_m\}$. Below we will also consider
the membership problem for 3-generator finite inverse monoids;
here the input consists of four injective partial maps $f_0, f_1,
f_2, f_3$, and the question is the same as before (now with $m =
3$).

{\sc PSpace}-completeness of the membership problem for general
functions was shown by D.~Kozen \cite{Kozen}. For permutations the
problem is in the complexity class NC (hence in P), as proved by
L.~Babai, E.~Luks, A.~Seress \cite{Luks}. In \cite{BMcTh}
M.~Beaudry, P.~McKenzie, D.~Th\'erien proved that the membership
problem for general functions (not assumed to be injective)
remains {\sc PSpace}-complete if the monoid generated by $\{ f_1,
\ldots, f_m\}$ is assumed to be in certain pseudo-varieties, and
is NP-complete or in NP or in P for certain other
pseudo-varieties.

Although inverse monoids are similar to groups in many ways,
problems about inverse monoids can be much harder than the
corresponding problems about groups:

\noindent
\begin{thm} \label{PSpace_compl_memb} \   
The membership problem for the class of finite {\rm inverse}
monoids is {\sc PSpace}{\em -complete.} The problem remains {\sc
PSpace}{\em -complete} if the finite inverse monoids are required
to have just three generators.
\end{thm}
{\bf Proof.} \ Since we showed that the intersection-emptiness
problem is {\sc PSpace}-complete for inverse finite automata with
a two-letter input alphabet, we can apply Kozen's reduction (see
p.~263 of \cite{Kozen}). Kozen's proof needs a few changes in
order to make his functions injective.

Let ${\cal A}_i = (Q_i, \Sigma, \delta_i, q_i^{\rm (start)},
q_i^{\rm (fin)})$ ($i = 1, \ldots, k$) be a sequence of inverse
finite automata, with the same two-letter alphabet $\Sigma =
\{\alpha, \beta\}$. We can assume that $q_i^{\rm (start)} \neq
q_i^{\rm (fin)}$ \ (see \cite{BMMW}). As the set acted on by our
partial functions we take
 \ $S = \{o_1, o_2\} \cup \bigcup_{i=1}^k Q_i$. The functions are defined
as follows:

For each $a \in \Sigma$, define $f_a : S \to S$ by $f_a(q_i) =
\delta_i(q_i, a)$ (for $q_i \in Q_i$), and $f_a(o_2) = o_2$.
However, $f_a(o_1)$ is undefined. Also, consider the function
$f_{\rm init} : S \to S$ defined by $f_{\rm init}(q_i^{\rm
(start)}) = q_i^{\rm (start)}$ for $i = 1, \ldots, k$, and $f_{\rm
init}(o_1) = o_2$, and $f_{\rm init}$ is undefined elsewhere.
Finally, the ``test function'' $f_0 : S \to S$ is defined by
$f_0(q_i^{\rm (start)}) = q_i^{\rm (fin)}$ for $i = 1, \ldots, k$,
and $f_0(o_1) = o_2$, and $f_0$ is undefined elsewhere.

Now it is straightforward to check (exactly as in \cite{Kozen},
p.~263) that $f_0$ is generated by \ $\{f_{\rm init}, f_{\alpha},
f_{\beta}\}^{\pm 1}$ \ iff \ \ $\bigcap_{i=1}^k L_{{\cal A}_i} \  
\neq \varnothing$.
 \ \ \ $\Box$


\bigskip

\bigskip

\noindent
{\bf Jean-Camille Birget} \\  
{\small Dept.\ of Computer Science, Rutgers University at Camden,
Camden NJ 08102, USA} \\  
{\tt birget@camden.rutgers.edu}

\smallskip

 and

\smallskip

\noindent
{\bf Stuart W.~Margolis} \\   
{\small Dept.\ of Mathematics, Bar Ilan University,
Ramat Gan 52900, Israel} \\  
{\tt margolis@macs.biu.ac.il}

\end{document}